\documentclass[preprint,5pt]{elsarticle}
\usepackage{multicol}
\usepackage{color}
\usepackage{xspace}
\usepackage{enumerate}
\usepackage{amssymb}
\usepackage{graphicx}
\usepackage{amsthm}
\usepackage{amsfonts}
\usepackage{amsmath}
\usepackage{MnSymbol}

\makeatletter
\def\bs{\expandafter\@gobble\string\\}
\def\lb{\expandafter\@gobble\string\{}
\def\rb{\expandafter\@gobble\string\}}
\def\@pdfauthor{Juan-Luis Domenech-Garret}
\def\@pdftitle{elsarticle.cls -- A documentation}
\def\@pdfsubject{Document formatting with elsarticle.cls}
\def\@pdfkeywords{LaTeX, Elsevier Ltd, document class}

\DeclareRobustCommand{\LaTeX}{L\kern-.26em%
        {\sbox\z@ T%
         \vbox to\ht\z@{\hbox{\check@mathfonts
           \fontsize\sf@size\z@
           \math@fontsfalse\selectfont
          A\,}%
         \vss}%
        }%
     \kern-.15em%
    \TeX}
\makeatother

\newdefinition{rmk}{Remark}
\newdefinition{example}{Example}

\newcommand{\vertiii}[1]{{\left\vert\kern-0.25ex\left\vert\kern-0.25ex\left\vert #1 
    \right\vert\kern-0.25ex\right\vert\kern-0.25ex\right\vert}}

\begin{document}

\def\testa{This is a specimen document. }
\def\testc{\testa\testa\testa\testa}
\def\testb{\testc\testc\testc\testc\testc}
\long\def\test{\testb\par\testb\par\testb\par}


\title{Numerical solution of several second-order ordinary differential equations containing logistic maps as nonlinear coefficients}

\author[1]{J.L. ~Domenech-Garret\corref{cor1}}\ead{domenech.garret@upm.es}
\author[2]{C.~Marin-Ferrer}\ead{carmen.marin@upm.es}
$•$
\cortext[cor1]{Corresponding author}

\address[1]{Dpto. F\'{\i}sica. 
E.T.S.I. Aeron\'autica y del Espacio.\\
Univ. Polit\'{e}cnica de Madrid,  Madrid, Spain.}
\address[2]{Dpto. Ing. Agroforestal. 
E.T.S.I.A.A.B.\\\
Univ. Polit\'{e}cnica de Madrid,  Madrid, Spain.}

\date{\today}
\begin{abstract}
This work is devoted to  find the numerical solutions of several one dimensional second-order ordinary differential equations. In a heuristic way, in such equations  the quadratic logistic maps regarded as a local function are inserted within  the nonlinear coefficient of the function as well as within the independent term. We apply the Numerov algorithm to solve these equations and we discuss the role of the initial conditions of the logistic maps in such solutions. 
\end{abstract}

\maketitle

\section{\textbf{Introduction}}
In this work, we will focus on the second order differential equations ($2^{nd}$-ODE) which contain the function, $f$, to be solved, their corresponding second derivative and an independent term. Such kind of equations can be written in a general form as

\begin{equation}
\label{2ODE}
f^{''}(x)\ +\  Q(x)\  f(x) = S(x)
\end{equation}
\noindent
where to shorten, as usual, we use $f^{'}, f^{''}, f^{'''}, ... ,$ to denote the first, second, third and higher derivatives of the function, $f$. Numerical methods applied to solve ODEs split the $x$ domain  into a lattice of  N points  according to
 a given step $h$. Such a discretisation, on practical grounds, transforms the original ODE into,
\begin{equation}
\label{2ODE-disc}
f_n^{''}\ +\  Q_n\  f_n = S_n
\end{equation}
\noindent
where, we write the  function $f_n\equiv f(x_{n})=f(x_{n-1}+h)$; $ Q_n\equiv Q(x_n)=Q(x_{n-1}+h)$;  and $ S_n\equiv S(x_n)=S(x_{n-1}+h)$. After discretisation, the algorithms to solve Eq.(\ref{2ODE-disc}) yield an appropriate recursion formula using educated guesses to get the initial values. Throughout this work, we will use the Numerov algorithm \cite{Numerov1,Numerov2,Koonin} to solve Eq.(\ref{2ODE-disc}). Such (forward) recursion formula, as we shall see, needs the knowledge of two initial values of $f$. Moreover, this formula also depends on the step, $h$, and the functional form of $Q_n$ and $S_n$. Briefly we write the recursion formula as,

\begin{equation}
\label{Rec1}
  f_{n+1}=F(f_n;f_{n-1}; h; Q_{n}; S_{n},...) \nonumber
\end{equation}

Here, as a main feature of this work, we take advantage of the numerical split of the domain to substitute $Q_n$ and $S_n$ by logistic maps, which themselves follow a recursive procedure \cite{Schuster, Nagashima}. Therefore, the functions $Q_n$ and $S_n$ correspond to functions of the form $A\ M(x)$, where $A$ is an appropriate constant magnitude, and $M(x)$ represents a quadratic logistic map. 
\begin{equation}
\label{Rec2}
  f_{n+1}=F(f_n;f_{n-1};h; M1(x_{n}); M2(x_{n}),...) \nonumber
\end{equation}
\noindent
The logistic maps, $M1, M2$, will then govern the diffusion of the magnitude $A$  along the domain. In such framework we set the correspondence between the recursive map and the ordinary discretisation of the functions $Q(x)$ and $S(x)$ when an ODE is solved numerically. As we will see later, the shape of such discretised $Q_n$ and $S_n$ along the entire domain $x_n$ is modelled by taking a particular logistic map, as well as their initial conditions. As we shall see, we find that the initial values of those maps actually define the equation to be solved.  At this point, we want to stress that in this effort we will not enter into rigorous aspects of the choice of such kind of coefficients. All these formal issues are beyond the scope of this numerical attempt, which is  performed  in a \emph{heuristic way}. This work is organized as follows: First, we review the Numerov method and the quadratic logistic maps to be used. In subsequent sections we provide the solution for several particular cases of $2^{nd}$-ODE  with different maps taking  $Q(x)=0$ and $S(x)=0$. Finally we will focus on the solution of the general case of the $2^{nd}$-ODE. Throughout this work, though  we take into account that $ Q_n $ refers to the non-linear coefficient and $ S_n $ to the independent term, we refer to both as coefficients.

\section{The Numerov method and the logistic maps  }
\label{NUM-M} 
The Numerov algorithm \cite{Numerov1,Numerov2,Koonin} provides a recursive formula to find the solution of the  $2^{nd}$-ODE. First, we split the $x$ range  into a lattice of  N points evenly spaced according to $x_{n}=x_{n-1}+h$ (where $h$ is the step). Then the $2^{nd}$-ODE step by step reads $f_n^{''}\ +\  Q_n\  f_n = S_n$. Subsequently, it is possible to expand the function around $x_n$ forward and backwards. The details about how to attain the Numerov formula, as well as the corresponding formula for the first derivative, can be found in the Appendix. Labelling $f_{n+1} \equiv f(x_{n}+ h)$ and $f_{n-1} \equiv f(x_{n}- h)$  we arrive to the following Numerov forward recursive relation, with a local error \textit{0}($h^6$):\\

\begin{eqnarray}
\label{NUMF}
f_{n+1}=\frac{2(1-\frac{5h^2}{12} Q_{n})\ f_{n}\ -\ (1+\frac{h^2}{12} 
Q_{n-1})\ f_{n-1} }{(1\ +\ \frac{h^2}{12} Q_{n+1})} + \nonumber \\ 
+ \frac{\frac{h^2}{12} \left(S_{n+1}+ 10 S_{n}+ S_{n-1}\right) }{(1\ +\ \frac{h^2}{12} Q_{n+1})} \: \: \: \: \: \: \: \: \: \: \: \: \: \: \: \: \: \: \: \: \: \: \: \: \: \: \: \: \: \: 
\end{eqnarray}
Following the same procedure, the corresponding backwards recursive relation can be attained. Therefore, to calculate the  solution function  using the recursive formulas also means to have knowledge of two initial values  for each direction. As in this effort we include logistic maps, we will need also to take into account their initial values to develop these recurrences. 

\begin{figure}[htb!]
\centerline{\hbox{ 
\includegraphics[width=\textwidth]{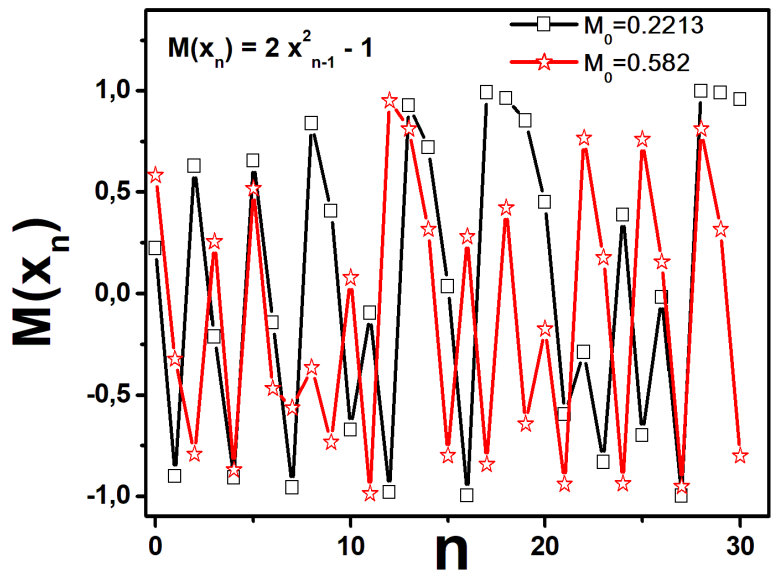} }}
\caption{(color on line): Map1 shapes along a given domain taking different initial values $M_0$.Line with squares $M_0= 0.2213$; Line with stars $M_0= 0.582$. }
\label{fig:Fig0}
\end{figure}

\subsection{\textbf{The logistic maps}}
Throughout this work we will take several quadratic logistic maps \cite{Schuster, Nagashima} as the nonlinear coefficients $Q_n$ and $S_n$ of the ODE. This will connect the different local values along the domain. These  are of the form $A\ M(x_{n})$, where $M(x_n)$ represents the recursive function that is dependent on the previous value through a given function $x_{n+1}= M(x_{n})$. Throughout this work, for the sake of simplicity, we normalize the constant $A$ to the unity. The selected quadratic maps we use here are labelled as  Map1 and Map2, given as:
\begin{itemize}
\item Map1
$$
x_{n+1}= k\ x^2_{n} + c 
$$
\noindent
where $k$ and $c$ are constants 

\item Map2
$$
x_{n+1}= B\  x_{n} \left( 1- x_{n} \right) 
$$
\noindent
where $B$  is a constant 

\end{itemize} 

The intrinsic nature of  the recursive maps means they need  the knowledge of the initial values to develop the recurrences for $Q_n$ and/or $S_n$. The solutions of the ODE (and, as we will see, the ODE itself) are governed by the form of the logistic map. Obviously, the  $Q_n$ and/or $S_n$ values need to be determined  along the whole domain prior to applying the Numerov formula given in Eq.(\ref{NUMF}). As mentioned above,  the solution for the function $f$ is strongly sensitive to the initial values of $Q_n$ and $S_n$ because these values determine the form of the nonlinear coefficients  along the $x_n$ domain. Figure \ref{fig:Fig0} displays an example of such behaviour: the same logistic Map1, taken here as $x_{n+1}= 2\ x^2_{n} -1$,  with different initial conditions evolves in a different way along the same domain.

\section{\textbf{Solutions of the discretised $2^{nd}$-ODE }}
 As  mentioned before, within this section we present the solutions of the function $f$ following different kinds of  the $2^{nd}$-ODE given by Eq.(\ref{2ODE-disc}), using the  logistic maps shown above. We start by looking at the particular case in which $Q_n=0$. Subsequently, we analise the one taking $S_n=0$. We finally find several solutions in the general case with non-zero $Q_n$ and $S_n$ coefficients. The Numerov recursion in each case  mentioned above is performed using the appropriate form of Eq.(\ref{NUMF}). All the results of this work have been computed using the  \textit{Python} language \footnote{The software used in this work is available upon reasonable request to the authors.}. A pedagogical approach of such a computer language can be found in \cite{PYTHON}. Those results have been as well checked using the PAW language \cite{PAW}.     
\begin{figure}[htb!]
\centerline{\hbox{ 
\includegraphics[width=\textwidth]{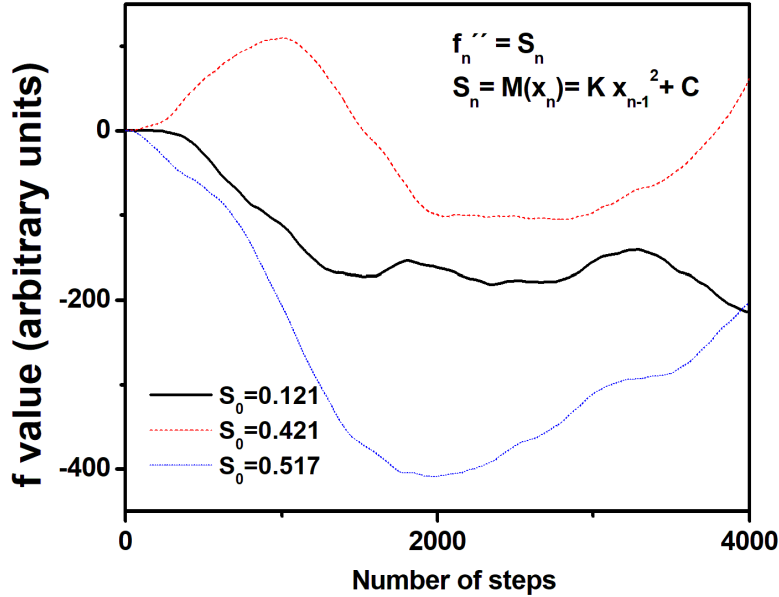} }}
\caption{(color on line): Solutions of the function f corresponding to the ODE $f_{n}^{''} = S_{n}$, where $S_n$ corresponds to a Map1 taking different initial values $S_0$. The $f$ solutions are: Upper line (dashed) taking $S_0= 0.421$; Middle (solid) line using $S_0= 0.121$; and lower (dotted) line with $S_0= 0.517$.  Here the step is $h=0.1$; $k=2, c=-1$. Initial values of $f$ to apply the Numerov algorithm are: $f_{0}=0.1; f_{1}=0.12$. }
\label{fig:Fig1}
\end{figure}

\subsection{ \textbf{Solutions of Equation $f_{n}^{''} = S_{n}$ ;\ ($Q_n=0$)} }
\label{Qeq0}
This kind of equations, in physics, can be cast in the framework of a given function $f(x)$ which relates it with its source $S(x)$. We can have a view of this case by taking the Poisson equation for the electric potential, $f$, and charge $S$ \cite{Koonin,Tijonov}. In our particular choice, the source (electric charge) could be viewed as obeying a recursive distribution along the space. By dropping all the 
$Q$ contributions, the Numerov formula for the function in this case is, 

\begin{equation}
\label{NumQ0}
f_{n+1}= 2 f_{n}\ - f_{n-1}\ +\ \frac{h^2}{12}\left(S_{n+1}+ 10 S_{n}+ S_{n-1}\right) 
\end{equation}

In Figure \ref{fig:Fig1}, the results of the solution $f$ of the Eq.(\ref{NumQ0}) using  a Map1 are shown. The graph displays three solutions for $f$ depending on the initial values, $S_0$, of the Map1. As mentioned above, we observe the dependence of the initial values (labelled as $S_0$) of the coefficient. $S_0$ will determine the form of the independent term, $S_n$, of the ODE along the complete domain.\\

This figure shows the result already mentioned: what is displayed is the solution of three differential equations. The coefficients of such equations  obey the same Map1, \mbox{$x_n=2x^2_{n-1}-1$}, but the initial condition changes in each case. There are actually different coefficients among them along the domain, and therefore this yields different equations to be solved. In each case the solution is a one-valued function.

\subsection{ \textbf{Solutions of Equation $f_{n}^{''} + Q_{n}\  f_{n} = 0$ ;\ ($S_{n}=0$):}}
\label{Seq0}
\begin{figure}[htb!]
\centerline{\hbox{ 
\includegraphics[width=8.0cm]{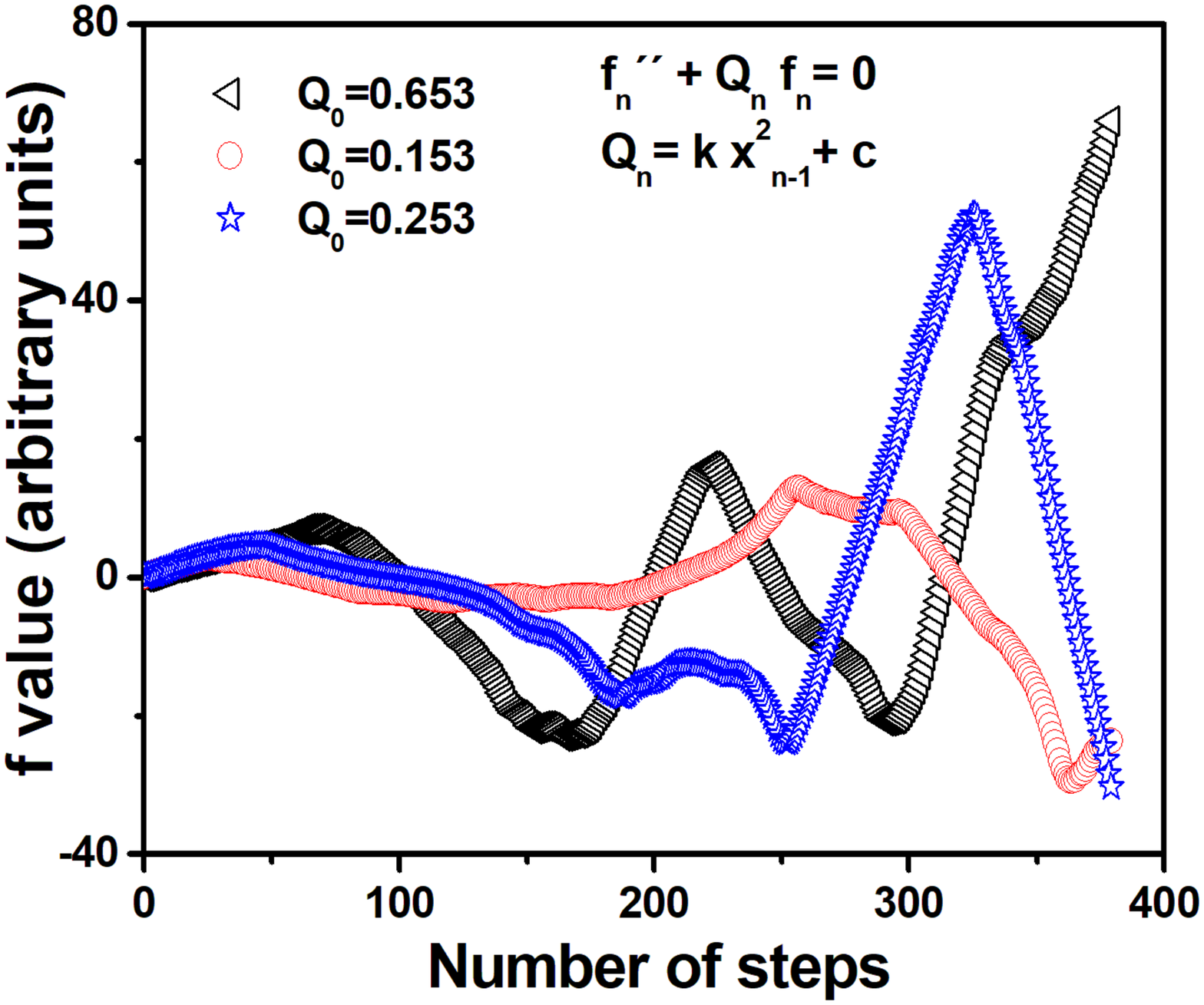} 
\includegraphics[width=8.0cm]{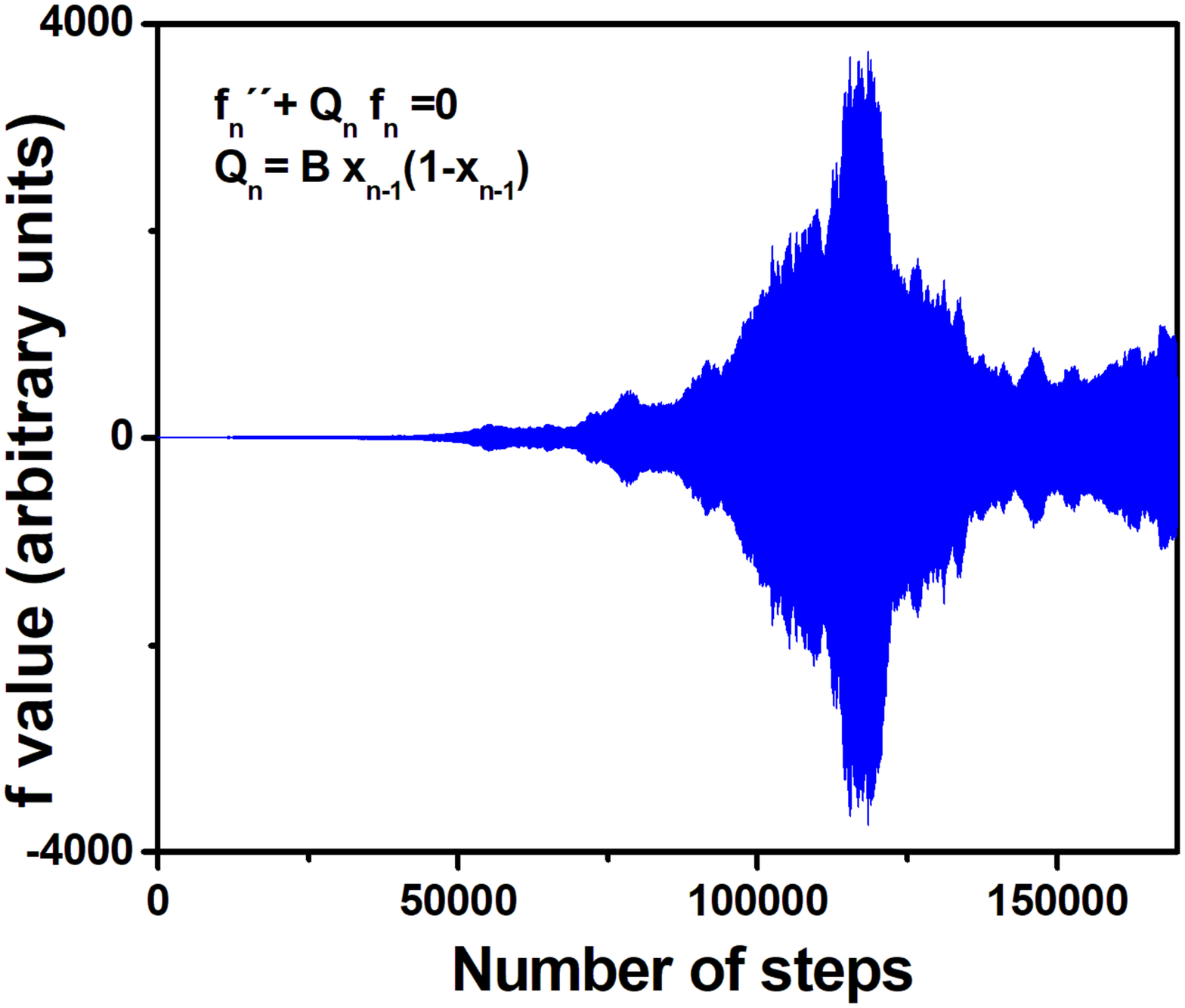}  }}
\caption{(color on line): Values of the function f corresponding to the ODE $f_{n}^{''} + Q_{n} f_{n} = 0 $. Left-(A): The  $Q_n$ coefficients correspond to the Map1 with $k=2,c=-1$. The step is $h=0.1$. Initial values: $f_0 = 0.1;\ f_1 = 0.2212$. Right-(B): $Q_n$ corresponds to the Map2 with $B=3.82$. Step $h=0.1$. Initial values: $f_0 = 0.1;\ f_1 = 0.1414;\ Q_0 = 0.93313$.}
\label{fig:Fig2}
\end{figure}
Equations having this form can be cast, for instance, in the framework of a Schrodinger-type equation \cite{Tijonov, Domenech}. Following such a picture, the potential is merged into the function $Q(x)$ \cite{Koonin}. Here, taking our choice for $Q$, the potential would follow a recursive distribution along the space. The Numerov formula for the function in this case, by dropping all of  the $S_n$ contribution, is, 

\begin{equation}
\label{NumS0}
f_{n+1}=\frac{2(1-\frac{5h^2}{12} Q_{n})\ f_{n}\ -\ (1+\frac{h^2}{12}
Q_{n-1})\ f_{n-1}}{(1\ +\ \frac{h^2}{12} Q_{n+1})}
\end{equation}

In Figure \ref{fig:Fig2}-A shows the solution  of three equations following Eq.(\ref{NumS0}) using the Map1 as  $Q_n$ with $k=2$ and $c=-1$. The relevant values inserted to obtain the solution are displayed within the caption. The initial values $Q_0$ yielding the respective equation to be solved can be viewed within the figure. In Figure \ref{fig:Fig2}-B we observe the results  of  Eq.(\ref{NumS0}) using  the Map2 within $Q_n$. The relevant values can be found in the caption. Again the shape of those particular $Q_n$ recursive maps, and hence the equations to be solved, are modelled by the corresponding  initial values of $Q_n$.

\begin{figure}[htb!]
\centerline{\hbox{ 
\includegraphics[width=8.0cm]{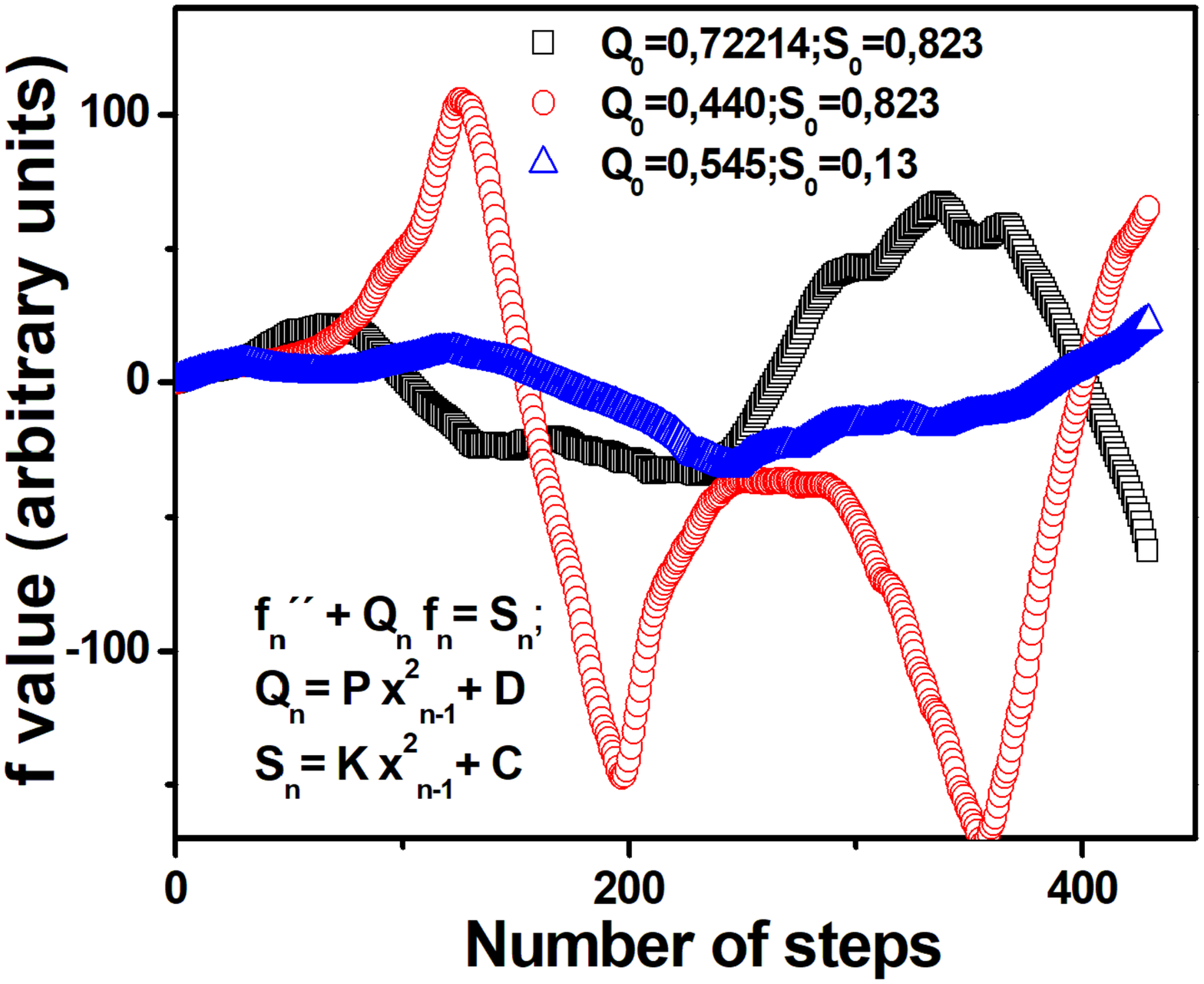} 
\includegraphics[width=8.0cm]{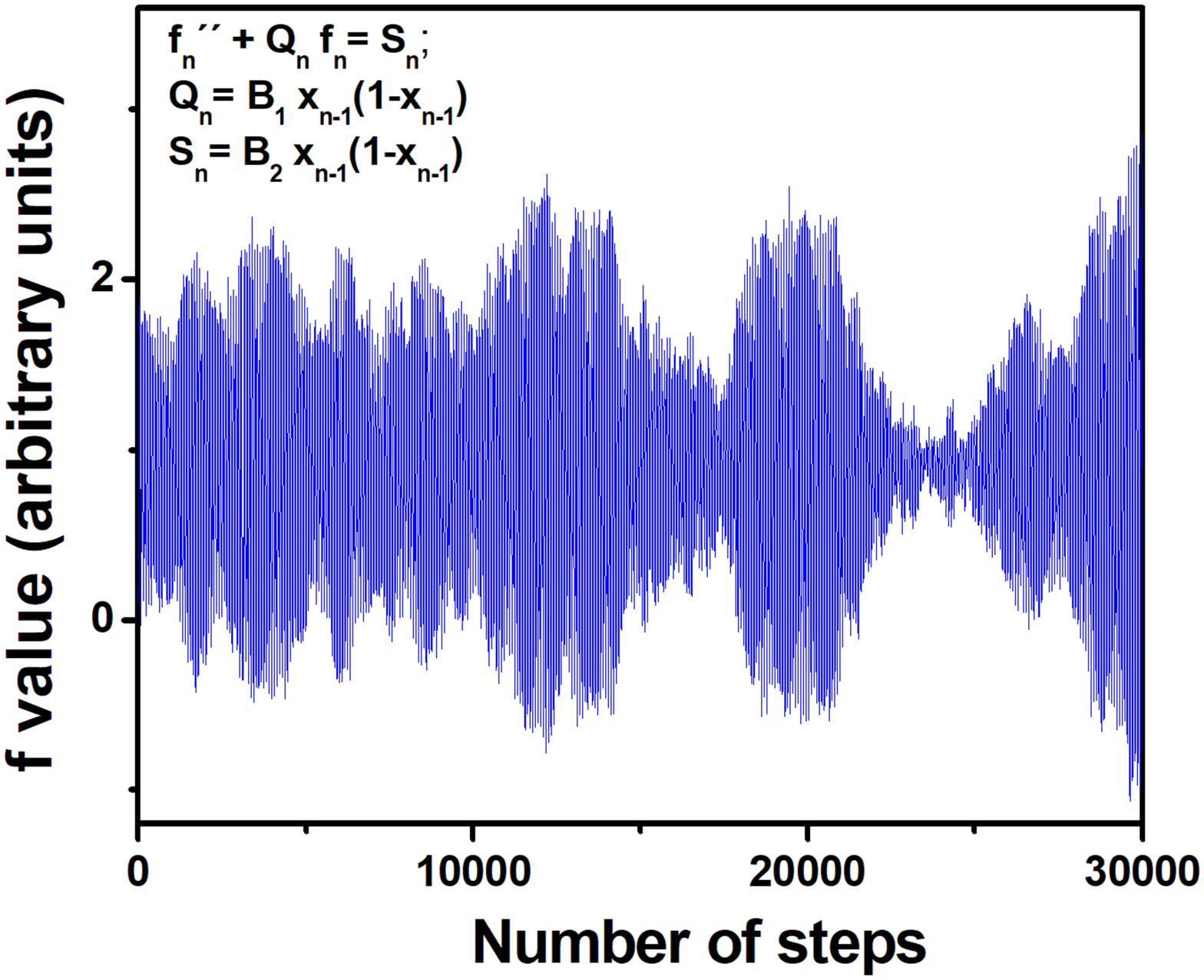} }}
\caption{(color on line):Values of the function f corresponding to the ODE $f_{n}^{''} + Q_{n} f_{n} = S_{n} $. Left (A): Both $Q_n$ and $S_n$ correspond to the Map1, both with $k=2,\ c=-1$. Step $h=0.1$. Initial values: $f_{0}=0.1;\ f_{1}=0.4212$.
Right (B): Both $Q_n$ and $S_n$ correspond to the Map2 with $B_1=3.342$ and $B_2=3.942$. Step $h=0.1$;  Initial values are: $f_{0}=0.1;\ f_{1}=0.1212;\ S_{0}=0.2821;\ 
Q_{0}=0.8213 $.}
\label{fig:Fig3}
\end{figure}

\subsection{ \textbf{Solutions of Equation $ f_{n}^{''} + Q_{n} f_{n} = S_{n} $. General case}}
\label{General}

This kind of equations can be found in several cases of interest in physics, for instance the elliptic equation  for the amplitude of permanent oscillations under the action of periodical forces \cite{Tijonov}. Figure \ref{fig:Fig3}-(A) shows the solution  of three equations following such a form. The $Q_n$ and $S_n$ coefficients obey  the Map1  with $k=2$ and $c=-1$. The relevant values inserted to obtain the solution can be found within the caption. The different initial values, $Q_0$ and $S_0$, yielding the respective equation to be solved are displayed within the figure. Again, the corresponding two initial values of $f_n$  required by the Numerov formula to find  such a solution can be found within the caption. In Figure \ref{fig:Fig3}-(B) we can observe the corresponding solution of an ODE taking the Map2 as coefficient both in $Q_n$ and $S_n$ using the initial values displayed within the caption of the Figure.\\

\begin{figure}[htb!]
\centerline{\hbox{ 
\includegraphics[width=8.0cm]{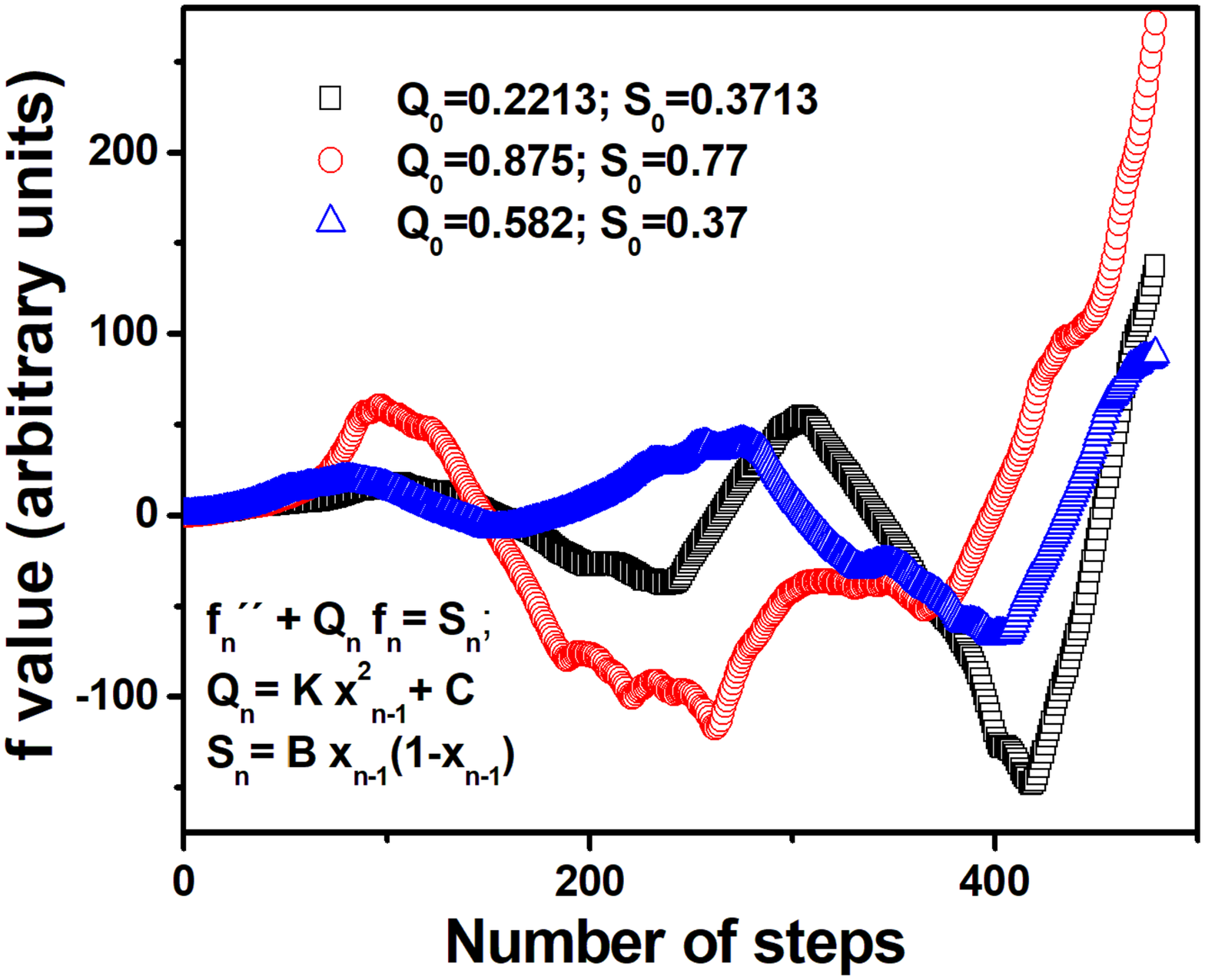} 
\includegraphics[width=8.0cm]{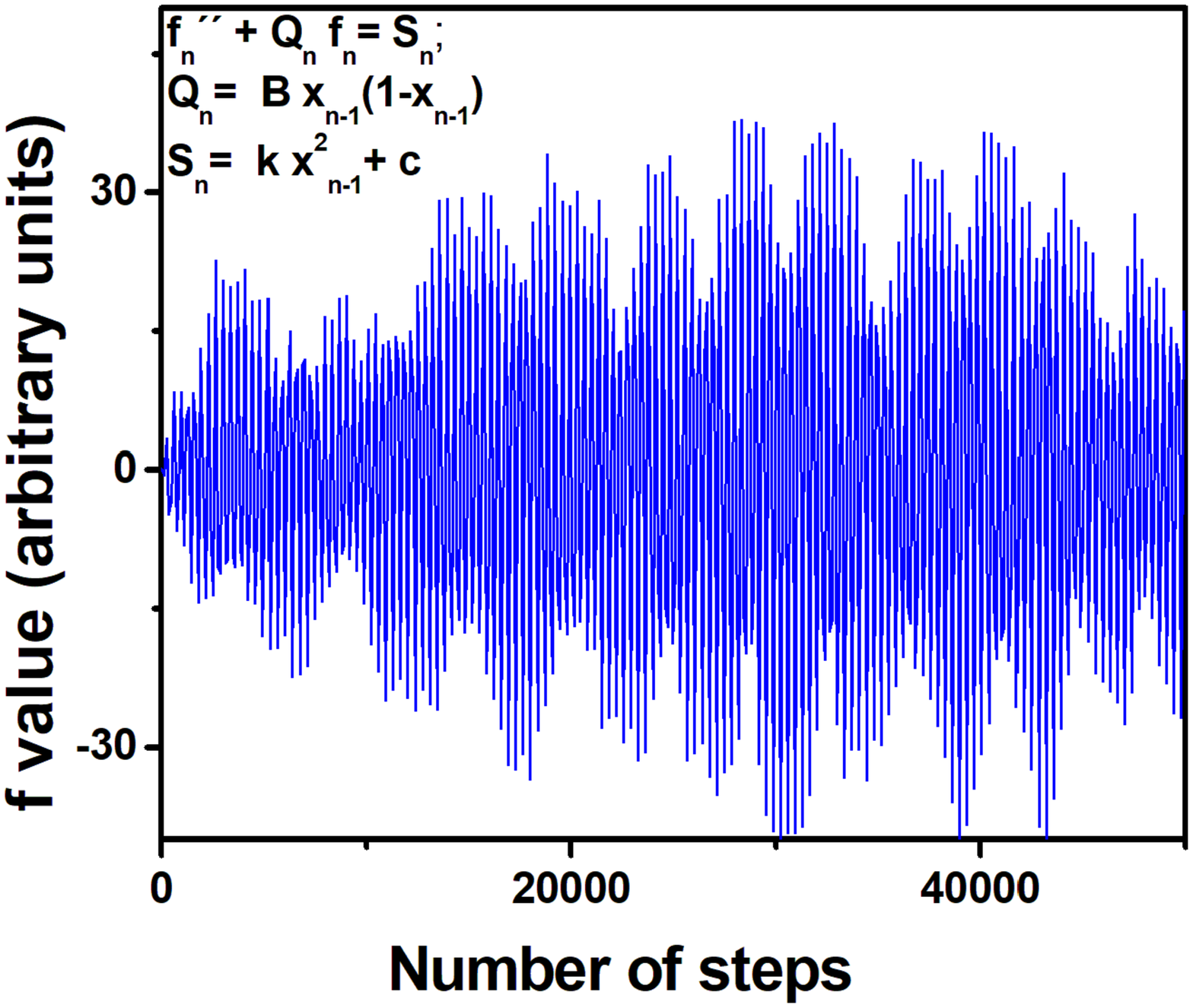}  }}
\caption{(color on line): Values of the function f corresponding to the ODE $f_{n}^{''} + Q_{n} f_{n} = S_{n} $. Left (A): $Q_n$  follows the Map1 with $k=2,c=-1$ and $S_n$ the Map2 with B=3.442. Step $h=0.1$. Initial values are: $f_{0}=0.1;\ f_{1}=0.15$. Right (B): $Q_n$  follows the Map2  with B=1.1342 and $S_n$ the Map1 with $k=2,c=-1$. Step $h=0.1$. Initial values are: $f_{0}=0.1;\ f_{1}=0.2212;\ Q_{0}=0.713;\ S_{0}=0.7213$.}
\label{fig:Fig4}
\end{figure}

 Figure \ref{fig:Fig4}-(A)  shows the solution taking  Map1 as the $Q_n$ coefficient and  Map2 within $S_n$. In Figure \ref{fig:Fig4}-(B) can be found the solution for $f$, in which $Q_n$ is a Map2, and the Map1 appears in $S_n$. In each case initial values $Q_0$ and $S_0$ de determine the shape of the  coefficients and therefore the ODE.  As well, each solution for $f$ is found  by taking the corresponding two initial values of the function, $f_0$ and $f_1$, needed using the Numerov algorithm.

\section{Conclusions.}
In this work we solve numerically several differential equations obeying the general form $ f_{n}^{''} + Q_{n} f_{n} = S_{n}$, in which a set of quadratic logistic maps are within nonlinear coefficients. In the general case, we also insert maps into the coefficients  having the same form but also taking different quadratic maps. Moreover, we find the solution of  equations taking $S_n=0$ and $Q_n=0$. In all of these cases, if we select a given Map and  set their  constants,  for instance the $k,c$ values in case of Map1, by taking a set of different initial values, we obtain the corresponding set of values for the coefficients along the domain. Because we deal with logistic maps, $M(x_n)$, the ODE above transforms into $f_{n}^{''} + M1(x_n) f_{n} = M2(x_n)$. Therefore, whenever we change the initial conditions of those non linear coefficients, we change the shape of those functions along the $x_n$ domain and actually  deal with a \emph{different} equation to be solved numerically.

\section{Acknowledgements} Authors acknowledge Prof. M. Marva-Ruiz for the useful comments, and Dr. J. Damba for the careful reading through the manuscript. J.L. ~Domenech-Garret acknowledges support by the Ministry of Science, Innovation and Universities of Spain under Grant number RT2018-094409-B-100.

\appendix
\section{\textbf{HOW TO ATTAIN THE NUMEROV FORMULA}}
Following section \ref{NUM-M}, up to $\textit{0}(h^5)$, the forward and backwards expansions of the function $f$ around $x_n$ are,

\begin{equation}
\label{fad}
 f_{n+1}= f(x_{n})+ h f'(x_n) +
\frac{h^2}{2}f^{''}(x_n) + \frac{h^3}{6}f^{'''}(x_n)+
\frac{h^4}{24}f^{(iv)}(x_n) 
\end{equation}

\begin{equation}
\label{fat}
f_{n-1} = f(x_{n})- h f'(x_n)+
\frac{h^2}{2}f^{''}(x_n) - \frac{h^3}{6}f^{'''}(x_n)+
\frac{h^4}{24}f^{(iv)}(x_n) 
\end{equation}
\noindent
By adding the two equations above, the odd terms cancel and then we obtain, up to $\textit{0}(h^6)$,

\begin{equation}
\label{sum}
f_{n+1}- 2f_{n}+ f_{n-1} = h^2 f_{n}^{''} +\ \frac{h^4}{12} f_n^{(iv)} 
\end{equation}
\noindent
Using Equation (\ref{2ODE-disc}), we can rewrite the last term of the above equation as follows:

\begin{equation}
\label{sum2}
f_{n}^{(iv)} = \frac{d^{2}f_{n}^{''}}{dx^2}= \frac{ d^{2} \left[S_{n} - Q_{n} f_{n}\right]  } {dx^2} \equiv S_{n}^{''}- Q_{n}^{''}\ f_{n}^{''}
\end{equation}
\noindent
If $F_n$ stands for a label for $Q_n, S_n, f_n,...$, approximating the second derivatives above by the three-point differentiation formula, \cite{3PDF1,3PDF2}
\begin{equation}
\label{sum3}
F_n^{''}\approx\frac{F_{n+1}- 2F_{n}+ F_{n-1}}{h^2} \nonumber 
\end{equation}
\noindent
And using Eq.(\ref{sum2}) into Eq.(\ref{sum}), after manipulations we finally attain the Numerov forward recursive relation, given by Eq.(\ref{NUMF}), with a local error \textit{0}($h^6$). 

\subsection{\textbf{Derivatives}} 

In order to have the corresponding formula of the first derivative at the appropriate order, we subtract  the expansions in Eqs. (\ref{fad}) and (\ref{fat}). We then get for the first derivative up to order \textit{0}($h^5$):
\begin{equation}
\label{1der}
f_{n}^{'}= \frac{1}{2h}\left[ f_{n+1} -  f_{n-1} -\ \frac{h^3}{3} f_n^{(''')}\right] 
\end{equation}
\noindent
On the other hand, using the above equation we can rewrite the third derivative as
\begin{equation}
\label{der3}
f_{n}^{'''} = \frac{d^{2}f_{n}^{'}}{dx^2}=\frac{1}{2h}\left[ f^{''}_{n+1} -  f^{''}_{n-1} -\textit{0}(h^3)\right]
\end{equation}
\noindent
and again, using the $2^{nd}$-ODE $f^{''}_{n}= S_{n} - Q_{n} f_{n}$ into the above equation, and including the result into Eq.(\ref{1der}) we finally attain, up to order $\textit{0}(h^4)$,

\begin{eqnarray}
\label{NumDer}
f^{'}_{n}= \frac{1}{2h} \left[ f_{n+1} -  f_{n-1}  \right] +  \nonumber \\
\frac{h^2}{6} \left[Q_{n+1}\ f_{n+1} - Q_{n-1}\ f_{n-1} - S_{n+1} + S_{n-1} \right] 
\end{eqnarray}
Therefore, with the knowledge of the two initial values and, again, the initial values $S_0$ and $Q_0$, since we include logistic maps,   the derivative can be obtained.


\begin{thebibliography}{99}

\bibitem{Numerov1} Fox L., Mayers D.F. (1987) \textit{Initial-value methods for boundary-value problems. In: Numerical Solution of Ordinary Differential Equations.} Springer, Dordrecht. https://doi.org/10.1007/978-94-009-3129-9-5.

\bibitem{Numerov2} J.W. Daniel and A.J. Martin (1977), \textit{Numerov's method with deferred corrections for two-point boundary-value problems}. SIAM J. Num. Anal., 14, 1033–105.

\bibitem{Koonin}    S.E.~ Koonin, D.C.~Meredith.  \emph{Computational physics}  Addison-Wesley, (1990) ISBN 0-201 -38623-2.

\bibitem{Schuster} H.G.~Schuster, \textit{Deterministic chaos: An Introduction}. VCH Verlagsgesellschaft mbH, Weinheim,
2nd edition, (1989).

\bibitem{Nagashima} H.~Nagashima and Y.~Baba \textit{Introduction to Chaos: Physics and Mathematics of Chaotic Phenomena.} IOP Publishing Ltd, (1999)  

\bibitem{PYTHON}  W.~McKinney \textit{Python for Data Analysis}  O'Reilly Media, Inc. Publisher ISBN: 9781449319793 (2012).

\bibitem{PAW}  \textit{PAW: Physics Analysis WorkStation} Avilable at CERN Program Library https://paw.web.cern.ch/paw/

\bibitem{Tijonov} A.N.~Tikhonov and A.A.~Samarskii \textit{Equations of mathematical physics } Ed. New York: Dover Publications, ISBN-13: 978-0486664224 (2011).


\bibitem{Domenech} J.L.~Domenech-Garret and M.A.~Sanchis-Lozano, Comput.\ Phys.\ Commun.\ {\bf 180},768
(2009). [arXiv:0805.2704 [hep-ph] ].

\bibitem{3PDF1} Chang Shu \textit{Differential Quadrature and Its Application in Engineering}, Springer, 2000, ISBN 978-1-85233-209-9.

\bibitem{3PDF2}  M.~Abramowitz and I. A. Stegun \textit{ Handbook of Mathematical Functions.} National Bureau of Standards, (1972).










\end{thebibliography}
\end{document}